\documentclass{amsart}

\usepackage[utf8]{inputenc}

\usepackage{amstext}
\usepackage{hyperref}
\usepackage{url}
\usepackage{caption}
\usepackage{a4wide}
\usepackage[utf8]{inputenc}
\usepackage{amsthm}
\usepackage{amsmath}
\usepackage{amsfonts}
\usepackage{amssymb}
\usepackage{natbib}
\usepackage{graphicx}
\usepackage{mathtools}
\usepackage{enumitem}
\usepackage{float}
\usepackage{placeins}
\usepackage{comment}
\usepackage{cutwin}

\newtheoremstyle{named}{}{}{\itshape}{}{\bfseries}{.}{.5em}{\thmnote{#1 }#1}
\theoremstyle{named}

\theoremstyle{plain}
\newtheorem{theorem}{Theorem}[section]
\newtheorem{proposition}[theorem]{Proposition}
\newtheorem{lemma}[theorem]{Lemma}
\newtheorem{corollary}[theorem]{Corollary}

\theoremstyle{definition}
\newtheorem{definition}[theorem]{Definition}

\newtheorem{example}[theorem]{Example}
\newtheorem{conjecture}[theorem]{Conjecture}

\newtheorem{remark}[theorem]{Remark}
\newtheorem*{theorem*}{Theorem}

\newenvironment{acknowledgments}{%
  \begin{abstract}
}{%
  \end{abstract}
}

\title[Explicit lower bounds on the conductor of abelian varieties with bad reduction]{Explicit lower bounds on the conductors of elliptic curves and abelian varieties over number fields}

\author{Pierre Tchamitchian}
\address{Institut de Mathématiques de Marseille, Université d’Aix-Marseille, 163, avenue de
Luminy, Case 907, 13288 Marseille Cedex 9, France}
\email{pierre.tchamitchian@univ-amu.fr}

\date{}


\begin{document}

\begin{abstract}

Following the work of Mestre, we use  Weil's explicit formulas to compute explicit lower bounds on the conductors of elliptic curves and abelian varieties over number fields with $L$-functions being``analytic $L$-functions".
Moreover, we obtain bounds for the conductor of elliptic curves and abelian varieties over $\mathbb{Q}$ with specified bad reduction and over number fields, when the same assumption is made on their $L$-function.
As an application, for specific fields, we prove the non-existence of abelian varieties, with analytic $L$-function, which have everywhere good reduction.
\end{abstract}

\maketitle


\section{Introduction}

A famous result of Fontaine~\cite{FontainevariétésabéliennesurZ} states that no abelian variety defined over $\mathbb{Q}$ can have good reduction everywhere. Mestre~\cite{Mestre} showed, under the assumption of standard conjectures that this result could be obtained again by using explicit formulas due to Weil~\cite{Weilformulesexplicites}, which relate exponential sums over primes with zeros of certain $L$-functions. More precisely this is done under the assumption that these $L$-functions satisfy certain conditions: they can be extended meromorphically to the complex plane, and satisfy a functional equation (precise conditions will be given in Section 2.1). Such function are said to be analytic $L$-functions. 
Applying these formulas to $L$-functions of abelian varieties, Mestre~\cite{Mestre} proved the following theorem:

\begin{theorem*}[Mestre, Proposition p.~21]\label{Theorem 1}
Let $F(x)$ be a test function (see Definition \ref{test function}) and let $\lambda >0$ be a real number. Then there is a computable bound $B_M(F,\lambda)$ such that for
any abelian variety over $\mathbb{Q}$ of dimension $g$ and conductor $N$, and whose $L$-function is an analytic $L$-function, we have \[ N \geqslant B_M(F,\lambda)^g \; .\] \end{theorem*}
For a particular choice of $\lambda$ and $F(x)$, Mestre obtained $B_M(F,\lambda) = 10.323$. In particular, this gives a conditional proof of Fontaine's result, under the standard conjecture that all abelian varieties over $\mathbb{Q}$ have analytic $L$-functions. This conjecture was proved in the case of elliptic curves as a corollary of Wiles' modularity theorem~\cite{Wilesmodularity}.

In this paper, we show that the formulas used by Mestre can be refined:

\begin{theorem*}[Theorem \ref{Mestre inequality abelian variety over Q}]
Let $F(x)$ be a test function (see Definition \ref{test function}) and let $\lambda >0$ be a real number. Let $N$ be positive integer. Then there is a computable bound $B(F,\lambda,N,g)$ such that for
any abelian variety over $\mathbb{Q}$ of dimension $g$ and conductor $N$, and whose $L$-function is an analytic $L$-function, we have \[ N \geqslant B(F,\lambda,N,g) \geqslant B_M(F,\lambda)^g.\] If $\lambda$ is sufficiently large, then  $B(F,\lambda,N,g) > B_M(F,\lambda)^g$.\end{theorem*}

The bound $B(F,\lambda,N,g)$ is determined using a finer analysis of the coefficients corresponding to primes of bad reduction of an abelian variety appearing in the Weil's explicit formulas. Consequently, we obtain better conductor bound than Mestre's for abelian varieties with analytic $L$-function with prescribed some primes of bad reduction. 

We also prove a generalization of a remark of Mestre~\cite[Remarque 3, p.~22]{Mestre}: in this remark Mestre explained that one obtains the non existence of elliptic curves with everywhere good reduction over certain quadratic fields. We extend this result to other number fields, and we obtain explicit examples for number fields of degree up to $12$. Moreover, this result generalizes to abelian varieties of higher dimension, as shown in the following theorem:

\begin{theorem*}[Theorem \ref{Principal th of the article}] Let $F(x)$ and $\lambda$ be as above. Let $K$ be a number field of degree $n$. There exists an explicitly computable real constant $B(K,F,\lambda)>0$, such that for any abelian variety defined over $K$ of dimension $g$ and conductor $\mathfrak{f}$, and with analytic $L$-function, the following inequality holds:
\[ \left( N_{K/\mathbb{Q}}(\mathfrak{f}) \right) ^{\frac{1}{n}} 
\geqslant B(K,F,\lambda)^g \; , \] where $ N_{K/\mathbb{Q}}$ denotes the norm of $K/\mathbb{Q}$.
\end{theorem*}

Using this result we obtain a list of $965$ number fields over which no abelian variety with analytic $L$-functions can have everywhere good reduction. The computation were done under the standard GRH conjecture. Without this hypothesis, we obtain a feasible weaker result: the number of fields goes down from $965$ to $545$.

Modularity theorems (Wiles proof of modularity of elliptic curves~\cite{Wilesmodularity}, Siksek and al~\cite{Freitasmodularity} work over real quadratic fields, etc) do provide sharp effective results in some particular cases, related to the computations of the relevant automorphic forms. Nevertheless, the results of this paper are applicable even if a modularity statement is unknown, and where reduction types are prescribed, specially over number fields.

This paper is organized as follows. In Section~2, we restate the explicit formulas due to Mestre in the formalism of analytic $L$-functions, as introduced by Farmer, Pitale, Ryan and Schmidt \cite{FarmerPitaleRyanSchmidt}. Using this formalism we then give Mestre's explicit formulas, closely following \cite{Mestre}.
In Section~\ref{Section refined formulas}, we give the explicit form of the constant appearing in theorem \ref{refined Mestre bounds} along with the proof of the Theorem, but only for elliptic curves. We then give some numerical examples of the improvements on the conductor bound given by the refined formulas. 

Section 4 analysis the form of the $L$-functions of abelian varieties, after which we derive conductor bounds for abelian varieties over $\mathbb{Q}$ as in Section 3.

Section 5 follows the structure of Sections 3 and 4, replacing $\mathbb{Q}$ with any number field. In particular, we investigate the relation between the constants $B(K,F,\lambda)$ and $B(L,F,\lambda)$ for any finite extension of number fields $L/K$. This is followed by numerical applications of the formulas and we provide data on number fields up to degree 16 that we found over which no abelian varieties with analytic $L$-function can have everywhere good reduction. These results are obtained under the Generalized Riemann Hypothesis (GRH). Unconditional results are possible, subject to a modification of the test function, giving weaker bounds.

An appendix where we discuss the choice of a test function is provided. In particular we present a method involving polynomial functions to improve on the bounds, exploiting the density of polynomials in the subspace of  test functions considered for the computations.
Implementation of the algorithms used to compute the conductor bounds and their results are publicly available on GitHub~\cite{DOI_conductor_abelian_varieties}.

\section{analytic \texorpdfstring{$L$}{L}-functions and Mestre's formulas}

We let $s=\sigma+it$ denote a complex variable, with $\sigma$ and $t$ real, and $f(s)$ to designate a function on the complex domain. Let $\Gamma(s)$ denote the usual Gamma function, and define the gamma factors $\Gamma_\mathbb{R}(s)$ and $\Gamma_\mathbb{C}(s)$ by 
  \[\Gamma_\mathbb{R}(s) = \pi^{-s/2}\Gamma(s/2)\qquad \mathrm{and} \qquad\Gamma_\mathbb{C}=(2\pi)^{-s}\Gamma(s) \cdot \]
If $f(s)$ is a complex function, its Schwartz reflection $\tilde{f}(s)$ is defined by $\tilde{f}(s)=\overline{f(\overline{s})}$.

\subsection{Analytic \texorpdfstring{$L$}{L}-functions}

Mestre's explicit formulas (Theorem \ref{Mestre's general formula}) involve $L$-functions which satisfy several specific conditions, which were highlighted by Mestre~\cite{Mestre}. The many similarities with the framework described in \citep[section 2.1]{FarmerPitaleRyanSchmidt} makes it interesting to use this formalism, in which the common properties of $L$-functions arising from diverse objects are described precisely. In particular, Mestre's explicit formulas can be applied to the $L$-functions arising from elliptic curves --- as done later in this paper --- but also modular forms, Galois representations, automorphic representations, and pure motives.  We recall two of their important properties:

Let $L(s)$ be an analytic $L$-function in the sense of \citep{FarmerPitaleRyanSchmidt}. There exists a positive integer $N$, called the \textbf{conductor}
of $L(s)$, a pair of non-negative integers $(d_1,d_2)$, and complex numbers $\{ \zeta_j\}^{d_1}_{j=1}$ and $\{ \nu_k\}^{d_2}_{k=1}$, such that the \textbf{completed $L$-function} 
\begin{flalign*}
      && \Lambda(s)= N^{s/2}\overset{d_1}{\underset{j=1}{\prod}}\Gamma_{\mathbb{R}}(s+\zeta_j)\overset{d_2}{\underset{k=1}{\prod}}\Gamma_{\mathbb{C}}(s+\nu_k)L(s) &&
\end{flalign*}
    satisfies the following functional equation: 
    \begin{flalign*}
      &&\Lambda(s)=\varepsilon \tilde{\Lambda}(1-s) \; , &&
\end{flalign*} for some  $\varepsilon \in \mathbb{C}$.

Furthermore, $L(s)$ admits a product expression, called the Euler product form of the $L$-function:

\begin{flalign*}
      && 
     L(s) = \underset{p \, \mathrm{prime}}{\prod}L_p({s})^{-1}= \underset{p \, \mathrm{prime}}{\prod}F_p(p^{-s})^{-1},
     &&
\end{flalign*}
absolutely convergent for $\sigma >1$, where the $F_p(T)$ are polynomials. The factors $L_p(s)$ are called the Euler factors.
Mestre's explicit formulas come essentially from the two previous points: the functional equation relates two sides of a complex integral on a closed path, and the Euler product form allows to explicitly compute one of the sides of the integral. From the residue theorem, we get a relation between a sum indexed over powers of primes and the zeros of an $L$-function. Other conditions in the definition of analytic $L$-functions are also needed for Mestre's explicit formulas to be valid, but are less intuitive than the two conditions highlighted above.

\subsection{Mestre's explicit formulas}
In the following, let $x$ denote a real variable. We first need to introduce the following definition:
\begin{definition}\label{test function} We say that $F: \mathbf{R} \rightarrow \mathbf{R}$ is a \emph{Mestre test function} if:
\begin{itemize}[noitemsep]

     \item[(i)] there exists $\varepsilon>0$ such that $F(x)\exp((1/2+\varepsilon)|x|)$ is integrable over $\mathbb{R}$,
 
     \item[(ii)] there exists $\varepsilon>0$ such that $F(x)\exp((1/2+\varepsilon)x)$ has bounded variation and is integrable over $\mathbb{R}$, and the value at each point is the average of its right and left limits,
  
     \item[(iii)] the function $(F(x)-F(0))/x$ has bounded variation over $\mathbb{R}$.
\end{itemize}

Any smooth function with compact support is a test function. Another example of test function is the Gaussian $F(x) = \exp(-x^2)$. 
The first condition  in their definition guarantees the existence of the associated Mellin transform\footnote{This is actually the Mellin transform of the function $v^{-{1}/{2}}F(\log v)$ defined for $v > 0$, see \cite[p.~4]{Weilformulesexplicites}.}

\begin{equation*}\label{Def Phi(s)}
\Phi_F(s) = \int_{-\infty}^{+\infty}F(x)e^{(s-1/2)x} \, dx
\end{equation*}

holomorphic on the vertical strip $\{s \in \mathbb{C},-\varepsilon  < \sigma < 1+\varepsilon \}$ where $\varepsilon$ is taken to satisfy conditions i) and ii). Such an $\varepsilon$ is said to be \textit{compatible} with $F(x)$. If $\varepsilon$ is compatible, then any $\varepsilon'$ with $0< \varepsilon' < \varepsilon$ is also compatible. The third condition is a technical conditions used for the proof of the formulas.
\end{definition}

\begin{theorem}[Mestre's explicit formulas {\cite[I.2]{Mestre}}]\label{Mestre's general formula} Let $F(x)$ be a Mestre test function and $\Phi(s)$ its associated Mellin transform. Let $L(s)$ be a analytic tempered $L$-function with completed $L$-function
\[ \Lambda(s) = N^{s/2}\overset{d_1}{\underset{j=1}{\prod}}\Gamma_{\mathbb{R}}(s+\zeta_j)\overset{d_2}{\underset{k=1}{\prod}}\Gamma_{\mathbb{C}}(s+\nu_k)L(s). \]
We define

\[ I(u,v) = \int_0^{\infty}  \big(  F(ux)e^{-(u/2+v)x}/(1-e^{-x})-F(0)e ^{-x}/x \big) \, dx \,,\]

\[ J(u,v) = \int_0^{\infty}\big( F(-ux)e^{-(u/2+v)x}/(1-e^{-x})-F(0) e^{-x}/x \big) \, dx  \,\cdot\]

\noindent We have, for $A = 2^{-d_2}N^{1/2} \pi^{-(d_1+2d_2)/2}$,
\[
\underset{\rho}{\sum}\Phi_F(\rho)
-\underset{\zeta}{\sum}\Phi_F(\zeta)
+\overset{d_1}{\underset{j=1}{\sum}}I(1/2,\zeta_j/2)+J(1/2,\overline{\zeta_j}/2)
+\overset{d_2}{\underset{k=1}{\sum}}I(1,\nu_k)+J(1,\overline{\nu_k})
\]
\[=2F(0)\log(A) - \underset{p,i,m \geqslant 1}{\sum}
(\alpha_{i,p}^mF(m\log p)+\overline{\alpha_{i,p}}^mF(-m\log p))
\frac{\log p}{p^{m/2}},\] where $\rho$ (resp. $\zeta$) runs over the set of zeros (rep. of poles) of $L(s)$ of real part between $0$ and $1$, counted with their multiplicities, for any $p$ the $\alpha_{i,p}$ are the roots of the polynomials $F_p(T)$ defining the local factor $L_p(s)$,  and the notation $\sum \Phi(\rho)$ means $\underset{T \rightarrow + \infty}{\lim} \; \underset{|\mathrm {Im}\rho | < T}{\sum}\Phi(\rho) $.
\end{theorem}

\begin{proof} The proof follows closely the presentation in \cite[I.2 Formules explicites]{Mestre}, with the following adaptations: 
\begin{itemize}
\item Condition (i) is given by Mestre without the absolute value since he considers only the application to even functions. Indeed, in this case, the absolute value is not needed.

\item The integral term  $I(u,v)$ differs from that in Mestre, who wrote \[I(u,v) = u\int_0^{\infty}  \big(  F(ux)e^{-(u/2+v)x}/(1-e^{-x})-F(0)e ^{-x}/x \big) \, dx \, .\] The first factor $u$ is a typo that does not affect the numerical applications, for which we have $u=1$ (both in Mestre~$\cite{Mestre}$ and this article).

\item The constant $A$ depends on the convention made for the normalization of $\Lambda(s)$ and its gamma factors.
\end{itemize}
\end{proof}

\subsection{Application of the formulas to elliptic curves}

Mestre used these formulas to prove that the conductor of an elliptic curve must be at least $11$. 
Let $E/ \mathbb{Q}$ be an elliptic curve defined on $\mathbb{Q}$ with conductor $N$ and associated $L$-function
  \[ L_E(s) = \underset{p \mid N}{\prod}(1-a_p p^{-s})^{-1}\underset{p \nmid N}{\prod}(1- a_p p^{-s}+p \cdot p^{-2s})^{-1} \; , \]
where $a_p = p+1-E(\mathbb{F}_p)$ if $p \nmid N$ and $a_p = 0$, $a_p = 1$ or $a_p = -1$ if $E$ has bad reduction at $p$ of additive type, split multiplicative type or non split multiplicative type, respectively. We define $\alpha_p$ and $\overline{\alpha_p}$ as the roots of the polynomial $F_p(T) = 1 - a_pT + p T^2$ when $p$ is a prime of good reduction. To apply the formulas, we need this function to be an analytic $L$-function:
\begin{proposition}\label{Prop modularity sur Q} The function $L(s)$ defined by\[ L(s) = L_E \left( s+ \frac{1}{2} \right)  \] is an analytic $L$-function.
\end{proposition}
\begin{proof} 
The functional equation is true for the $L$-functions of certain modular forms~\cite[Theorem~5.10.2]{Diamond-Shurman} and Wiles' modularity theorem~\cite[Theorem~0.4, p.~448]{Wilesmodularity} identify these $L$-functions with those of elliptic curves. Here the change of variables ensures the functional equation is correctly normalized. 

The bound on the coefficients comes from the Hasse-Weil bound when $p$ is prime of good reduction and is a straightforward computation when $p$ is a bad prime (see \cite[Theorem~V.2.3.1 and Appendix~C.16]{SilvermanI}). 
\end{proof}
When $L(s)$ is as above, its completed $L$-function is 
\[ \Lambda(s) = N^{s/2}\Gamma_\mathbb{C}\left(s+\frac{1}{2}\right)L(s). \] 
Finally, we assume that $L(s)$ satisfies the Generalized Riemann Hypothesis (GRH) and the Birsch Swinnerton-Dyer conjecture (BSD). Explicitly, the zeros of $\Lambda(s)$ are assumed to lie on the line $ \{ s \in \mathbb{C}, \, \Re(s) = \frac{1}{2} \} $, and the order at $s = \frac{1}{2}$ of $L(s)$ is the rank $r$ of $E$.

\medskip
To get Mestre's result, we need stronger hypotheses on the test functions we use.
\begin{definition}\label{compact test function definition}
We say that $F(x)$ is a \emph{compact Mestre test function} if $F(x)$ is an even Mestre test function, with support in $[-1,1]$, positive Fourier transform, and with $ F(0)=1$. 
\end{definition}

Since $F(x)$ has positive Fourier transform, by Fourier reciprocity we obtain $F(0)> 0$. Therefore the choice $F(0) = 1$ leads to no loss of generality.
Such functions are for example obtained by considering a function $g(x)$, even, $\mathcal{C}^1$ compactly supported on $[-\frac{1}{2},\frac{1}{2}]$, such that $\int_\mathbb{R}g(x)^2\,dx = 1$. Then, setting $F(x) = g * g(x)$ gives us a compact Mestre test function.

In the following, for any real number $\lambda >0$ and compact Mestre test function $F(x)$, we define $F_\lambda(x) = F(\frac{x}{\lambda})$. Such a function is essentially a compact Mestre test function, but with support in $[-\lambda,\lambda]$. In particular its Fourier transform is $ \widehat{F_\lambda}(x) = \lambda \hat{F}(\lambda x)$.

\begin{proposition}\label{fundamental inequality}
Let $F(x)$ be a compact Mestre test function, and $\lambda >0$. Let $r$ be the rank of $E$. Then,

\[
\lambda r \Phi_{F_\lambda}(\frac{1}{2})  + 2\underset{p,m \geqslant 1}{\sum}
b_m(p)\frac{F_\lambda(m\log p)}{p^m}\log p
   + M_{\lambda,F} \leqslant \log(N) \; ,
   \]
where 
\begin{equation*}\label{Def M_lambda}
M_{\lambda,F} 
= 2\log(2\pi)+ 2I_{F_{\lambda},2} 
= 2 \left( \log(2\pi)+ \int_0^{\infty} \frac{F_\lambda(x)e^{-x}}{1-e ^{-x}}-\frac{e^{-x}}{x} \, dx \right) \cdot 
\end{equation*}

\end{proposition}
\begin{proof}
If $F(x)$ is a compact Mestre test function, then the explicit formulas applied to $L(s)$ give 
\[
\underset{\rho}{\sum}\Phi_{F_\lambda}(\rho)
+ 2\underset{p,m \geqslant 1}{\sum}
b_m(p) \frac{F_\lambda(m\log p)}{p^m}\log p
=\log(N)-2\log(2\pi)- 2I(1,1/2) \; , \]
where $b_m(p)=a_p^m$ if $p \mid N$ and $b_m(p) = \alpha_p^m + {\overline{\alpha_p}}^m$ otherwise. We remark that since $F(x)$ is even, the integral terms $I(u,v)$ and $J(u,v)$ in Theorem \ref{Mestre's general formula} are equal.
We will note
\[ 
I_{F} = I(1,1/2) = \int_0^{\infty} \frac{F(x)e^{-x}}{1-e ^{-x}}-\frac{e^{-x}}{x} \, dx  \; \cdot 
\]
\noindent Therefore, we can rewrite the last equality as
\[
\underset{\rho}{\sum}\Phi_F(\rho)+ 2\underset{p,m \geqslant 1}{\sum}b_m(p)\frac{F(m\log p)}{p^{m}}\log p =\log(N)-2\log(2\pi)- 2I_{F}   \; \cdot 
\]
By compacity of the support of $F(x)$ we have  
\[ \underset{p,m \geqslant 1}{\sum}
b_m(p)\frac{F_\lambda(m\log p)}{p^m}\log p = \underset{p^m \leqslant e^\lambda}{\sum}
b_m(p)\frac{F_\lambda(m\log p)}{p^m}\log p \, \cdot\]
We isolate the zero of imaginary part $t=0$, which is of order $r$ assuming BSD, obtaining the term $\lambda r \Phi_{F_\lambda}(\frac{1}{2})$.
Finally, by the positivity of $\hat{F}(t)$ the Fourier transform of $F(x)$ we get \[ \underset{\rho \neq \frac{1}{2}}{\sum}\Phi_F(\rho) \geqslant 0 \, , \] since 
\[ 
\Phi_{F_\lambda} \left( \frac{1}{2}+it \right) 
= \lambda \hat{F}(\lambda t) \geqslant 0\; .
\] Note that here we used GRH for the general form of the zeros of $L(s)$. 
\end{proof}

\begin{remark}\label{without GRH}
This choice of a positive Fourier transform in definition \ref{compact test function definition} was used together with GRH in the proof of Proposition \ref{fundamental inequality}. We indicate now how to obtain a similar proposition without supposing GRH on the $L$-functions we work with. The point is to choose a test function $F(x)$ such that \[ \mathrm{Re}(\underset{\rho \neq \frac{1}{2}}{\sum}\Phi_F(\rho)) \geqslant 0 \, , \] which we then leads to an almost identical proof of Proposition \ref{fundamental inequality}. Without GRH, the only information we have on the $\rho$'s is that they are contained in the critical vertical strip $\{ s = \sigma +it \in \mathbb{C}, \;  \frac{-1}{2} < \sigma < \frac{1}{2} \}$. To ensure the inequality, one can simply makes sure that $\mathrm{Re}(\Phi_F(s))$ is positive for all $s$ in the critical band. This is achieved by the following trick, attributed to Odlyzko in~\cite[p.~8]{Poitou2}: take $F(x)$ a Mestre test function with positive Fourier transform. Then the function $G(x) = F(x)/\cosh(x/2)$ is also a test function such that $\Phi_G(s)$ satisfy the desired property.
\end{remark}

From this theorem, Mestre used the Weil bound to produce a conductor lower bound which only depends on the rank of the elliptic curve $E$.
\begin{corollary}
[{\cite[p.~8]{Mestre}}]\label{Mestre theorem}Let $E / \mathbb{Q}$ be an elliptic curve of rank $r$ and conductor $N$. For any $\lambda>0$ and $F$ a compact Mestre test function, the following inequality holds: 
\[\lambda r \Phi_F(\frac{1}{2}) - 2\underset{p^m \leqslant e^\lambda}{\sum}
\frac{ \left\lfloor 2\sqrt{p^m}  \right\rfloor \left| F_\lambda (m\log p) \right| }{p^{m}} 
\log p
+ M_{\lambda,F}
\leqslant \log(N)   \cdot \] 
%
\end{corollary}

\begin{proof} If $p | N$ then $b_m(p)$ belongs to $\{-1,0,1 \} $, and if $p \nmid N$ one has $b_m(p) = \alpha_p^m + \overline{\alpha_p}^m$. Since $\alpha_p^m$ and $\overline{\alpha_p}^m$ are both of absolute value $p^{m/2}$, we obtain $|b_m(p)| \leqslant 2p^{m/2}$ for all $m \geqslant 1$. Moreover $b_m(p)$ is always an integer, hence 
\[ b_m(p) \geqslant - \left\lfloor 2\sqrt{p^m}  \right\rfloor  \cdot\]
 The result follows by Proposition \ref{fundamental inequality}.
\end{proof}
Note that the only terms that depends on $E$ on the LHS is $r$. This allows one to get a lower bound on the conductor of any elliptic curve with conductor $r$. 

\section{The refined explicit formulas for elliptic curves over $\mathbb{Q}$}\label{Section refined formulas}

\subsection{The refined formulas}
 
In the proof of Corollary~\ref{Mestre theorem}, the Weil bound which was used on the coefficients $b_m(p)$ were independent of whether or not $p$ was a prime of bad reduction. We can give better estimates by separating the cases of bad and good reduction. In the following, we assume that every considered analytic $L$-function satisfy GRH.


\begin{theorem}\label{refined Mestre bounds}Let $E / \mathbb{Q}$ be an elliptic curve of rank $r$ and conductor $N$. Let $S_1$ denote the set of primes dividing $N$ exactly, $S_2$ that of primes whose square divides $N$, and finally $S = S_1 \cup S_2$. For any $\lambda>0$ and $F$ a positive compact Mestre test function, the following inequality holds: 

\[ \lambda r \Phi_F(\frac{1}{2})
- 2 \underset{\substack{p \notin S \\ p \leqslant e^\lambda}}{\sum}
\left(
\frac{  \left\lfloor 2 \sqrt{p^m} \right\rfloor F_\lambda(m\log p)}{p^{m}}\log p 
\right)
- 2 \underset{\substack{p \in  S_1 \\ p \leqslant e^\lambda}}{\sum}
 \left(
\frac{F_\lambda(m \log p)}{p^m}\log p
\right)
 +M_{\lambda,F} 
 \leqslant
    \log(N)  \cdot \]

\end{theorem}

\begin{proof}
  The left hand side of the inequality in Proposition \ref{fundamental inequality}
\[
\lambda r \Phi_F(\frac{1}{2})  + 2\underset{p,m \geqslant 1}{\sum}
b_m(p)\frac{F_\lambda(m\log p)}{p^m}\log p
   + M_{\lambda,F} \leqslant \log(N) \]

\noindent can be separated into three sums depending on whether $p \in S_1$, $p\in S_2$ or $p \notin S$. If $p$ does not divide $N$, then the inequality \[ b_m(p) \geqslant - \left\lfloor 2\sqrt{p^m}  \right\rfloor  \] is optimal, but if $p|N$ then we simply have $ b_m(p) \geqslant -1 $, and similarly if $p^2|N$ then $ b_m(p) = 0 \, ,  $ hence the result.
\end{proof}

This theorem holds under the same assumption as in \cite{Mestre}. Nevertheless, the result gives improved bounds when the reduction types are specified. 
\subsection{Numerical applications for elliptic curves}

We now illustrate the computational improvements obtained with Theorem~\ref{refined Mestre bounds}. To do so we compare them with the extensive data we have on elliptic curves over $\mathbb{Q}$, available in the LMFDB~\cite{lmfdb}. Although no new result appears, this will provide a measure of how precise this formula can be, and how one can expect sharp conductor bounds in higher dimensions cases.

Let $F(x)$ be the following Mestre test function: \[ F(x) = (1-|x|) \cos(\pi x) + \frac{\sin(\pi x)}{\pi} \quad \textrm{on $[-1,1]$, and zero elsewhere.} \] 
Hereafter we refer to this function as the \emph{Odlyzko test function}. Since this function is twice derivable and with compact support, it is automatically a test function. Its Fourier transform is positive, as the Odlyzko funciton is the convolution square of the function in $t$ defined $\cos(\pi t)$ on $[-\frac{1}{2},\frac{1}{2}]$ and zero elsewhere.
Using the Odlyzko test function, as Mestre did,  we present some computations in Tables~\ref{table of rank 0 curves} and \ref{table rank  1 curves} (in particular, these bounds are valid under GRH). The tables are structured as follow:

\begin{description}
    \item[First column] Displays primes we suppose to be of good reduction.
    
    \item[Second and third columns] The second column lists primes of multiplicative reduction, while the third column gives primes of additive reduction.
    
    \item[Fourth column] Provides the best choice of a parameter $\lambda$, to within ~$0.01$, for obtain the optimal bound.
    
    \item[Archimedean bound] The next entry gives the \textbf{archimedean bound} $B_\mathbb{R}$, computed using the test function and the corresponding $\lambda$, followed by its computation time.
    
    \item[Bound $B_\mathbb{Z}$] The second last column contains the bound $B_\mathbb{Z}$, obtained from $B_\mathbb{R}$ and the congruence condition on the conductor:
    \begin{itemize}
        \item A prime of good reduction does not divide the conductor.
        \item A prime $p$ of multiplicative reduction divides the conductor once.
        \item A prime $p$ of additive reduction divides the conductor with exponent $2$, except when $p=2$ or $p=3$, where the exponent can be larger.
    \end{itemize}
    
	\item[Exact bound] The lowest conductor of an existing elliptic curve with the required conditions. Such a conductor is found, for example, by searching through the LMFDB~\cite{lmfdb} for curves or rank lower than or equal to $5$. 
	
\end{description}

We note that the more we give primes of bad reduction, the less the bound $B_\mathbb{R}$ is relevant, as the congruence conditions become much stronger. For this reason, we only suppose there are few small primes of bad reduction for the examples in the tables.	

\begin{table}
\begin{tabular}{|c|c|c|c|c|c|c|c|} \hline
good  & mult & add & $\lambda$ & $B_\mathbb{R}$ & comp. time & $B_\mathbb{Z}$ & exact bound      \\ \hline \hline

$3$      &     $2$   &        &    $1.47$ & $12.956$ & 0.011820 s & 14  & 14  \\ \hline
 
$2$      &    $3$    &        &    $1.49$ & $10.823$ & 0.007873 s& 15  & 15    \\ \hline

$2$, $3$ &      $5$  &        &    $1.31$ & $10.394$ & 0.008945 s& 35  & 35  \\ \hline
 
$2$, $3$ &   $7$     &        &    $1.31$ & $10.394$ & 0.006919 s& 35  & 35   \\ \hline

$3$      &           & $2$    &    $1.68$ & $17.293$ & 0.009318 s& 20  & 20  \\ \hline

$2$      &           &  $3$   &    $1.69$ & $11.552$ & 0.006477 s& 27  & 27    \\ \hline

$2$, $3$ &           & $5$    &    $1.31$ & $10.394$ & 0.007274 s& 25  & 25    \\ \hline

	     &  $2$, $3$ &        &    $1.80$ & $15.037$ & 0.007443 s& 30  & 30    \\ \hline

$3$      &  $2$, $5$ &        &    $1.47$ & $12.956$ & 0.010621 s& 70  & 70    \\ \hline
 
         &     $2$   &    $3$ &    $1.98$ & $17.444$ & 0.006577 s& 18  & 54    \\ \hline

$3$      &     $2$   &   $5$  &    $1.47$ & $12.956$ & 0.006659 s& 50  & 50    \\ \hline

         &    $3$    & $2$    &    $1.99$ & $22.525$ & 0.006668 s& 24  & 24    \\ \hline
 
$3$      &     $5$   & $2$    &    $1.69$ & $17.302$ & 0.005707 s& 20  & 20    \\ \hline
  
\end{tabular}
    \caption{Bounds for rank 0 curves}
    \label{table of rank 0 curves}
\bigskip
\begin{tabular}{|c|c|c|c|c|c|c|c|} \hline

good 		& mult 		& add 	& $\lambda$ & $B_\mathbb{R}$ & time & $B_\mathbb{Z}$ & exact bound  \\ \hline \hline

			&  			    &  			& $1.68$ & $34.566$ & 0.014996 s& 35  	&	37		\\ \hline

 $3$		&      2     	&  			& $1.96$ & $51.571$ & 0.007215 s& 58    	&	58		\\ \hline

$2$, $5$	&      3     	&  			& $1.99$ & $44.926$ & 0.006920 s& 51 	& 	57   	\\ \hline

$2$, $3$	&      5     	&  			& $1.69$ & $34.588$ & 0.007134 s& 35		&	65    	\\ \hline
	
$2$, $3$	&      7      	&  			& $1.68$ & $34.566$ & 0.026686 s& 35		&	77    	\\ \hline

$3$			&       		&  2    	& $2.15$ & $81.113$ & 0.014097 s& 88		&  	88    	\\ \hline

$2$			&       		&  3   		& $2.15$ & $54.927$ & 0.010600 s& 63		&  	99   	\\ \hline
	
$2$, $3$	&     			&  5        & $1.70$ & $34.598$ & 0.006389 s& 175	&	175    	\\ \hline
	
$2$, $3$	&      			&  7        & $1.68$ & $34.566$ & 0.006172 s& 49		&	245    	\\ \hline
   
			&      2, 3     &			& $2.30$ & $79.323$ & 0.012507 s& 102	&	102   	 \\ \hline
	
$3$			&      2, 5     &   		& $2.24$ & $55.34$  & 0.017157 s& 70  	&	130	  	\\ \hline

$2$			&      3, 5     &   		& $2.27$ & $48.76$  & 0.011524 s& 75		&	135    	\\ \hline

			&      2 		& 3   	 	& $2.47$ & $104.971$ & 0.011536 s& 126 	&	162   	\\ \hline

$3$			&      2 		& 5    		& $2.39$ & $58.434$ & 0.012172 	s& 350 	&	350   \\ \hline

			&      3 		& 2    		& $2.49$ & $137.242$ & 0.010039 s& 156	&	156    	\\ \hline

$2$			&      3 		& 5    		& $2.45$ & $51.797$ & 0.011728 	s& 75	&	525    \\ \hline

$3$			&      5 		& 2    		& $2.48$ & $95.06$ & 0.010690  s& 140	&	160    	\\ \hline

$2$			&      5 		& 3    		& $2.53$ & $64.737$ & 0.011297 s& 135  	&	135  \\ \hline

$2$, $3$	&      5 		& 7    		& $1.69$ & $34.588$ & 0.007112 s& 245	&	245    \\ \hline
   
			&      			&   2,3 	& $2.65$ & $189.709$ & 0.010291 s& 216	&	216		 \\ \hline

\end{tabular}
    \caption{Bounds for rank 1 curves}
    \label{table rank  1 curves}    
\end{table}

\bigskip

We now give some remarks related to the effectiveness of the refined formula, while discussing some of the results in the tables.
\begin{itemize}

    \item
    Table~\ref{table of rank 0 curves} presents results for rank $0$ curves. The minimal conductor for such curves is $11$, but assuming bad reduction at small primes significantly improves the lower bound. Assuming multiplicative (resp. additive) reduction at $2$ gives a lower bound of $12.956$ (resp $17.293$). Combined with congruence conditions, this means such curves have conductor greater than 14 (resp 20), and indeed curves with these conductor do exist. In particular, there are no curves of conductor $16$. 

    \item
    For rank $1$ curves (Table~\ref{table rank 1 curves}) we obtain that with no constraints, the lower bound is $35$, very close to the actual minimal conductor $37$. Assuming for example additive reduction at $2$ and good reduction at $3$ gives an archimedean lower bound of $81$, which becomes $88$ after applying congruence conditions. This example, with the previous paragraph, shows how it is possible to obtain much higher bounds which are still close to the reality.

    \item
    The method becomes more effective as the assumed rank increases. For rank $1$, improvements from the explicit formula are more pronounced, and the congruence conditions play a smaller role. 

    \item
    The congruence conditions are often strong enough on their own, even without improvements from the explicit formula. For example, assuming additive reduction at $5$ and multiplicative reduction at $3$ gives $N \equiv 0 \pmod{75}$, which leads to a much stronger lower bound than Mestre’s method ($N \geqslant 11$). Supposing bad reduction at primes larger than $3$ (e.g., $5$ or $7$) produces negligible improvements, as only small primes $p < e^\lambda$ affect the exponential sum in the explicit formulas.

    \item
    There are failures to obtain an optimal lower bound. For example this occurs when assuming additive reduction at $5$ for one rank $1$: the explicit formula suggests that a curve of conductor $N = 25$ could exist, but in fact, no such curve does. The same phenomenon can be observed when assuming multiplicative reduction at $2$ and additive reduction at $3$ for rank $0$ curves. For primes $p > 3$, assuming bad reduction has little influence, though one can notice a slight increase in the optimal parameter $\lambda$. Overall, the method confirms that only the smallest primes have a significant impact on the lower bounds obtained.
    
    \item Finally, we remark that the optimal parameter $\lambda$ seems to increase with the choice of a higher rank and more pronounced bad reduction. For abelian varieties of higher dimension, which tend to have higher rank than elliptic curves, the choice of an optimal parameter $\lambda$ is therefore expected to be relatively large. In particular, more small primes will have an influence on the conductor bound (as all primes smaller than $e^\lambda$ appear in the explicit formula).

\end{itemize}

\newpage
\section{The refined formulas for abelian varieties over $\mathbb{Q}$}\label{section ab var over Q}

We now indicate how similar improvements can be made when studying abelian varieties. Let $A$ be an abelian variety of dimension $g$ and conductor $N$, defined over $\mathbb{Q}$. We associate to $A$ its first étale cohomology group and the correspondent $L$-function $L_A(s)$ that conjecturally is an analytic $L$-functions, up to an affine change of variable: see \citep[5.6, in particular 5.6.14]{Kahn}. 

Therefore we apply Theorem \ref{Mestre's general formula} on $L_A(s)$. Nevertheless, to compute results on $A$ similar to the case of elliptic curves, we first give a precise expression of the local factors of $L_A(s)$, depending on the type of bad reduction of $A$ at each prime. Mestre used a similar approach in \cite[p.~19]{Mestre} without giving details, and we found this result hard to extract from the existing literature. Hence, we state it and prove it below.
\subsection{The $L$-function of an abelian variety} 

Let us note $\mathcal{A}_p$ the identity component of the special fiber at $p$ of  the Néron model $\mathcal{A}$ of $A$. By Chevalley's theorem there is an exact sequence
\[ 0 \rightarrow U \times T \rightarrow  \mathcal{A}_p \rightarrow B \rightarrow 0 \]  where $B$ is an abelian variety over $\mathbb{F}_p$, $T$ is an algebraic torus and $U$ is a unipotent group, of respective dimensions $g_{ab}(p)$, $g_m(p)$ and $g_u(p)$ such that $g_{ab}(p) + g_m(p)+g_u(p) = g$. Saying that $A$ has good reduction at $p$ means that $\mathcal{A}_p$ is an abelian variety, or equivalently that $g_{ab}(p)=g$ and $g_m(p)=g_u(p)=0$.

\begin{proposition}\label{form of L function} The $L$-function $L_A(s)$ associated to the first étale cohomological group of $A$ is an Euler product with each local factor $L_p(s)$ of the following form: when $p \nmid N$,
\[ L_p(s) = \underset{j = 1}{\overset{g}{\prod}}((1- \alpha_j p^{-s})(1- \overline{\alpha}_j p^{-s}))^{-1}  \]
where $|\alpha_j| = \sqrt{p}$, and when $p|N$,
\[ L_p(s) = \underset{j = 1}{\overset{g_{ab}}{\prod}}\big((1- \alpha_j p^{-s})(1- \overline{\alpha}_j p^{-s})\big)^{-1}
\prod_{j=1}^r \Phi_{m_j} ( p^{-s})^{-1} \; \; ,\]
where $|\alpha_j| = \sqrt{p}$ and $\Phi_m$ is the $m$-th cyclotomic polynomial. In the latter case we have $g_m(p) = {\sum_j}\varphi(m_j)$ and $g_{ab}(p) + g_m(p) \leqslant g$ .

\end{proposition}
For $g=1$ we recover the usual form of the local Euler factors of elliptic curves defined over $\mathbb{Q}$.  
\begin{proof}

We follow the method indicated in \cite[\S 2.3]{Serre_facteurs_locaux}: to compute the $L$-function of an abelian variety is to compute the characteristic polynomial of the action of the geometric Frobenius on its étale cohomology. 
More precisely, fix distinct primes  $p$ and $\ell$. We consider the scalar extension of $A$ to $\mathbb{Q}_p$, still denoted by $A$. Let $\overline{\mathbb{Q}}_p$ be a separable closure of $\mathbb{Q}_p$, $\mathrm{G}_p = \mathrm{Gal} ( \overline{\mathbb{Q}}_p / \mathbb{Q}_p)$, $\mathrm{Frob}$ the arithmetic Frobenius and I the inertia subgroup of G.  Set
\[ V =  H^1(A_{\overline{\mathbb{Q}}_p},\mathbb{Q}_\ell)^{\mathrm{I}} \, , \] 
the subgroup fixed by I. 
The local factor $L_p(s)$ is defined by 
\[ L_p(s) = \det(1-p^{-s}\mathrm{Frob}^{-1} | H^1(A_{\overline{\mathbb{Q}}_p},\mathbb{Q}_\ell)^{\mathrm{I}}) .\]

We use several isomorphisms to make the computations explicit. First, we have the following isomorphisms of G-modules using \cite[Theorem~15.1]{Milne1986} and the Weil pairing
\[ H^1(A_{\overline{\mathbb{Q}}_p},\mathbb{Z}_\ell) \cong \mathrm{Hom}(T_{\ell}(A),\mathbb{Z}_{\ell}) \cong T_\ell(A^\vee)(-1)  \]
where $A^\vee$ denotes the dual variety of $A$,
\[ T_\ell(A^\vee)(-1) = T_\ell(A^\vee) \otimes \mathbb{Z}_\ell(-1)  \] 
and $\mathbb{Z}_\ell(-1) = T_\ell(\mathbb{G}_m)^*$. There is a non-canonical, Galois equivariant isomorphism of $\mathbb{Q}_\ell$-modules between $V_\ell(A)$ and $V_\ell(A^\vee)$ given by any polarization of $A$. Therefore, the characteristic polynomial of action of Frob$^{-1}$ on $H^1(A_{\overline{\mathbb{Q}}_p},\mathbb{Z}_\ell)$ is the same as the one of its action on $V_\ell(A)(-1)$.

The action of  Frob on  $\mathbb{Z}_\ell(1)$ corresponds to multiplication by $p$, as does the action of the geometric Frobenius on $\mathbb{Z}_\ell(-1)$: by definition $T_\ell(\mathbb{G}_m)$ is the projective limit of the $\ell^n$ roots of unity  $\varprojlim \mu_{\ell^n}$ in $\overline{\mathbb{Q}}_p$ and Frob acts by raising its elements to the $p$-power. Dualising and taking the inverse of Frob, we obtain the same result for the action of the geometric Frobenius on $\mathbb{Z}_\ell(-1)$. Thus, for any vector space $V$ with an action of Frob, the characteristic polynomial of the action Frob on $V(-1)$ is that of the action of $p.$Frob on $V$, by definition of the tensor product of representations.

The inertia subgroup acts trivially on both $\mathbb{Z}_\ell(1)$ and $\mathbb{Z}_\ell(-1)$ (because its action on $l^n$ roots of unity is trivial), so, tensoring with $\mathbb{Q}_\ell$, we obtain
\[  H^1(A_,\mathbb{Q}_\ell)^\mathrm{I} \cong V_\ell(A)^\mathrm{I}(-1).\]
Recall that $\mathcal{A}_p$ is the identity component of the special fiber of the Neron model of $A$. By \cite[Lemma $2$]{SerreTate} and the finiteness of the component group, we have
 \[T_\ell(A)^\mathrm{I} \cong T_\ell(\mathcal{A}_p) \, \cdot\] 

By Chevalley's theorem we have the exact sequence
\[ 0 \rightarrow U \times T \rightarrow  \mathcal{A}_p \rightarrow B \rightarrow 0 \] over $\mathbb{F}_p$, where $U$ is unipotent, $T$ is a torus and $B$ is an abelian variety. Therefore, taking the Tate-modules, one finds
\[ 0 \rightarrow T_\ell(U) \times T_\ell(T) \rightarrow  T_\ell(\mathcal{A}_p) \rightarrow T_\ell(B) \rightarrow 0. \] since multiplication by $\ell$ in  $U \times T$ is surjective over $\overline{\mathbb{F}}_p$, as $\ell \neq p$. Note that this sequence is automatically exact on the left and in the middle by the definition of the Tate module, and only the exactness on the right needs $U\times T$ to be $\ell$-divisible.
Thus we only need to compute the characteristic polynomial of the action of the geometric Frobenius on the three Tate modules $T_\ell(B)$, $T_\ell(U)$ and $T_\ell(T)$.

The $L$-factor of an abelian variety over $\mathbb{F}_p$ is well known, see \citep[Proposition 2.1.2]{SerreCurveFiniteFields}, so it remains to treat the case of the unipotent group and the torus.

The case of the unipotent group is also simple: since $\ell \neq p$, there are no $\ell$-torsion points in $U$, and so $T_\ell(U) \cong \{ 0\}$. The corresponding $L$-factor in $L_A(s)$ is trivial.

We are left with the study of the action of Frobenius on the Tate module of a torus over $\mathbb{F}_p$. 
Let $\mathbb{F}/\mathbb{F}_p$ be a finite extension of degree $n$ such that that $T$ is isomorphic to some $\mathbb{G}_{m,\mathbb{F}}^d$ over $\mathbb{F}$. Hence,
\[ T_\ell(T) \cong T_\ell(\mathbb{G}_{m,K}^d)\] as Gal$(\overline{\mathbb{F}}_p/\mathbb{F})$ modules. In particular, the action of $\mathrm{Frob}^n$ is that of multiplication by $p^n$ for the $\mathbb{Z}_\ell$-module structure of $T_\ell(\mathbb{G}_{m,\mathbb{F}}^d)$, and therefore, the action of $(p \, \mathrm{Frob}^{-1})^n$ is trivial. So $\mathrm{Frob}^{-n}$ acts trivially on $T_\ell(T)(-1)$, and the characteristic polynomial of $\mathrm{Frob}^{-1}$ on $V_\ell(T)(-1)$ must be a product of cyclotomic polynomials, which concludes the proof.

\end{proof}

\subsection{The refined formulas}
The completed $L$-function of $A$ is defined by 
\[ \Lambda(s) = N^{s/2} \Gamma_{\mathbb{C}}(s)^gL_A(s) \]
and conjecturally verifies 
\[ \Lambda(s) = \pm \Lambda(2-s) \cdot\]


This, together with Proposition \ref{form of L function}, immediately leads to the following theorem: 
\begin{theorem}\label{Mestre inequality abelian variety over Q}
Assume that $L(s)= L_A(s+ \frac{1}{2})$ is an analytic $L$-function satisfying GRH. Let $S$ the set of primes dividing $N$. Let $F(x)$ be a compact Mestre test function and $\lambda>0$ a real number. Let $r$ denote the analytic rank of the abelian variety $A$. Then  
\[ \lambda r \Phi_F(\frac{1}{2})
- 2 g \underset{\substack{p \notin S \\ p^m \leqslant e^\lambda}}{\sum}
\left(
\frac{  \left\lfloor 2 \sqrt{p^m} \right\rfloor F_\lambda(m\log p)}{p^{m}}\log p 
\right)
- 2  \underset{\substack{p \in  S \\ p^m \leqslant e^\lambda}}{\sum}
 g_m(p) \left(
\frac{F_\lambda(m \log p)}{p^m}\log p
\right)
 \]
 \[ - 2  \underset{\substack{p \in  S \\ p^m \leqslant e^\lambda}}{\sum}
g_{ab}(p) \left(
\frac{  \left\lfloor 2 \sqrt{p^m} \right\rfloor F_\lambda(m\log p)}{p^{m}}\log p 
\right)
 + g M_{\lambda,F}
 \leqslant
  \log(N)  \cdot \]
\end{theorem}

As in the section for elliptic curves, a slight change on the test functions we use in the theorem lead to a similar result which does not GRH, see remark~\ref{without GRH}.

\subsection{Numerical applications for abelian varieties}

Specifying the reduction type at certain primes leads again to more precise conductor bounds. the  We say that bad reduction at $p$ is of type $(k,n,m)$ when $g_{ab}(p) =k$, $g_m(p) = n$ and $g_u(p)=m$. In the tables below, we display conductor bounds for abelian surfaces with reduction of type $(0,n_2,m_2)$ at $2$ and reduction of type $(1,n_3,m_3)$ at $3$, for rank $0$ and rank $1$ abelian surfaces respectively. The column $B_\mathbb{R}$ gives the computed bound.

\begin{table}[H]
\begin{tabular}{|c|c|c|c|} \hline
type at $p=2$ 		& type at $p=3$ 			& $\lambda$ & $B_\mathbb{R}$ \\ \hline \hline

	(0,2,0)		&  	(1,1,0)	    				& $1.62$ & 185.9 	\\ \hline

 	(0,2,0)		&   (1,0,1)    				& $1.71$ & 202.4 	\\ \hline

	(0,1,1)		&   (1,1,0)  	  			& $1.73$ & 258.2 	\\ \hline

	(0,1,1)		&   (1,0,1)        			& $1.82$ & 289.8 	\\ \hline
	
	(0,0,2)		&   (1,1,0)        			& $1.83$ & 371.7 	\\ \hline

	(0,0,2)		&   (1,0,1)    				& $1.92$ & 428.9  	\\ \hline

\end{tabular}
	\caption{Inequalities for abelian surfaces of rank $0$ with bad reduction at $2$ and~$3$}
    \label{table abelian varieties of rank 0}    
\end{table}

\begin{table}[H]
\begin{tabular}{|c|c|c|c|} \hline
type at $p=2$ 		& type at $p=3$ 			& $\lambda$ & $B_\mathbb{R}$ \\ \hline \hline

	(0,2,0)		&  	(1,1,0)	    				& $1.87$ & 768.4 	\\ \hline

 	(0,2,0)		&   (1,0,1)    				& $1.97$ & 899.3 	\\ \hline

	(0,1,1)		&   (1,1,0)  	  			& $1.98$ & 1160.8	\\ \hline

	(0,1,1)		&   (1,0,1)        			& $2.06$ & 1398.3	\\ \hline
	
	(0,0,2)		&   (1,1,0)        			& $2.07$ & 1810.0 	\\ \hline

	(0,0,2)		&   (1,0,1)    				& $2.15$ & 2232.4 	\\ \hline

\end{tabular}
	\caption{Inequalities for abelian surfaces of rank $1$ with bad reduction at $2$ and~$3$}
    \label{table abelian varieties of rank 1}    
\end{table}

Supposing bad reduction at certain primes, without specifying the type of bad reduction, still allows us to give relevant bounds. Indeed, the weakest bound among all types of bad reduction is obtained when $g_{ab}(p) = g-1$, $g_m(p) = 1$ and $g_a(p)=0$, but is already better than the one given by Mestre. This allows us to give the following table for example. 
The second column displays the bound given by Mestre, which is $10.323^g$.

\begin{table}[h]
\begin{tabular}{|c|c|c|c|c|} \hline
dimension 	&  Mestre's estimate &   p = 3   & p = 5 		    &   p = 7	  \\ \hline \hline

2		    &  108.0			&	140.0	    &  	132.4			&  	132.4 	\\ \hline

3			&  1086.3		&	1426.1 	    &	1370.4			&	1370.4		  \\ \hline

4			&  11168	.2		&	14687.3 		&	14205.7		    &	14205.7		  \\ \hline

\end{tabular}
    \caption{Bounds for abelian varieties with bad reduction at $2$ and $p$}
    \label{Bounds for abelian varieties}
\end{table}

These results were obtained using the Odlyzko test function, with an optimal parameter $\lambda$ found in a discrete interval of values. Again, these bounds using these results are valid under GRH.

\begin{remark} In this section we used an inequality on the coefficient $b_m(p)$ to produce explicit formulas that are relevant for every abelian varieties. When working with one specific abelian variety, trying to have as much information on the form of the local factors of its $L$-factors as possible, and using Theorem~\ref{Mestre's general formula} yields more accurate conductor bounds. Indeed, the term $\underset{\rho \neq 1/2}{\sum}\Phi_{F_\lambda}(\rho)$ in Theorem~\ref{Mestre's general formula} tends to $0$ as $\lambda$ goes to infinity. Every other terms appearing in the formulas is computable, under the condition that one knows exactly what are the local factor for all $p \geqslant e^\lambda$, and what is the analytic rank of the abelian variety.

\end{remark}

\section{Application to abelian varieties over number fields}

We now turn our focus on elliptic curves and abelian varieties over number fields. While the question of conductor bounds in the case of prescribed bad reduction can be treated just as in the previous section, we focus only on the existence of elliptic curves and abelian varieties with everywhere good reduction.

Until the end of the section, we assume that every analytic $L$-function we consider satisfy GRH. Remark \ref{without GRH} will easily adapt to complement the theorems below is we don't assume GRH.

\subsection{Elliptic curves over number fields}
Let $K$ be a number field and $E/K$ an elliptic curve. The $L$-function attached to $E$ is defined as 
\[ L_{E}(s) = \prod L_{\mathfrak{p}}(s)^{-1} \] where, if we write $q = N_{K/\mathbb{Q}}(\mathfrak{p})$,

\begin{itemize}
\item[$\bullet$]$L_{\mathfrak{p}}(s) = 1 - a_{\mathfrak{p}}q^{s} + q^{1-2s}$ if $E$ has good reduction at $\mathfrak{p}$,
\item[$\bullet$]$L_{\mathfrak{p}}(s) = 1 - a_{\mathfrak{p}}q^{s}$ if $E$ has bad reduction at $\mathfrak{p}$, where in this case the value of $a_\mathfrak{p}$ becomes

\begin{align*}
a_{\mathfrak{p}} = 
\left\{
    \begin {aligned}
          +1 & \quad\textrm{if $E$ has split multiplicative reduction at $\mathfrak{p}$}, \\
          -1 & \quad\textrm{if $E$ has nonsplit  multiplicative reduction at $\mathfrak{p}$}, \\
           \  0 \ & \quad \textrm{if $E$ has additive reduction at $\mathfrak{p}$}.                  
    \end{aligned}
\right.
\end{align*}
\end{itemize} 

Set $L(s) = L_E \left( s+\frac{1}{2} \right)$.
The associated completed $L$-function is 
\[ \Lambda(s) = c^{s/2}\Gamma_\mathbb{C}(s)^n L (s) \; ,\] where $n=\deg K$ and $c$ is defined by 
\[ c = N_{K/\mathbb{Q}}(\mathfrak{f}_{E/K})d_K^2\] where $\mathfrak{f}_{E/K}$ is the conductor of the curve $E$, $N_{K/\mathbb{Q}}$ the norm map  of $K$ and $d_K$ the absolute discriminant of $K$. 

\medskip
We assume the following conjecture, known as the Hasse-Weil conjecture, for the rest of the paper.

\begin{conjecture}[Hasse-Weil] The function $L(s)$ satisfy the functional equation in the definition of an analytic $L$-function.
\end{conjecture}

This conjecture has been proved in some cases: when $K = \mathbb{Q}$ it follows from the modularity theorem, proved first by Wiles \cite{Wilesmodularity} for semi-stable elliptic curves. More recent works established the modularity of elliptic curves over $K$ in the case where $K$ is a real quadratic field (see \cite{Freitasmodularity}), a totally real cubic field (see \cite{Derickxmodularitycubicfields}) or totally real quartic field not containing $\sqrt{5}$ (see \cite{box2021ellipticcurvestotallyreal}). In our context, the Hasse-Weil conjecture provides the analytic continuation of the function $L_{E}(s)$ and the functional equation
\[L_{E}(2-s)= wL_{E}(s) \] with $w =\pm1$. The other axioms defining an analytic $L$-function can be proved using the definition of $L_{E}(s)$ and the Weil bound, as for the case $K=\mathbb{Q}$.
The explicit formulas should then be applicable to the $L$-function of any elliptic curve defined over a number field, leading to a generalization of the formula for elliptic curves over $\mathbb{Q}$. 

\begin{proposition}\label{equality in number fields}
Let $K$ be a number field of degree $n$ and discriminant $d_K$, $E$ a curve defined over $K$ of conductor $f_{E/K}$, $L_{E}(s)$ the $L$-function of $E$. Let $\alpha_\mathfrak{p}$ and $\alpha'_\mathfrak{p}$ be the reciprocal roots of $L_\mathfrak{p}$ when $E$ has good reduction at $\mathfrak{p}$. Let $F(x)$ be a Mestre test function with $F(0) =1$, and let finally $\lambda >0$. Set $b_m(\mathfrak{p})=a_\mathfrak{p}^m$ if $p \mid \mathfrak{f}_E$ and $b_m(p) = \alpha_\mathfrak{p}^m + {\alpha'_\mathfrak{p}}^m$ otherwise. Then the following holds:
\[ \lambda\underset{\rho}{\sum}\Phi_F(\rho) + 2\underset{\mathfrak{p},m \geqslant 1}{\sum}b_m(\mathfrak{p})\frac{F_\lambda(m\log q)}{q^{m}}\log q+ n M_{\lambda,F} =\log(N_{K/\mathbb{Q}}(\mathfrak{f}_{E/K})d_K^2),  \]
where $\rho$ runs over the set of zeros of the function $L_E(s+\frac{1}{2})$, $\mathfrak{p}$ runs over the prime ideals of $K$ and $q=N_{K/\mathbb{Q}}(\mathfrak{p})$, and $M_{\lambda,F}$ as in Proposition \ref{fundamental inequality}

\end{proposition}

\begin{proof} While this is formally the same formula as in the case $K = \mathbb{Q}$, we provide an alternative proof. We show that $L(s)$ is an analytic $L$-function. The functional equation is already assumed, and the Euler form of $L(s)$ is given by definition: for a fixed $p$, the product $L_p(s)$ of the factors $L_\mathfrak{p}(s)$ for all $\mathfrak{p}$ above $p$ has the desired form given in the definition of analytic $L$-functions. 
Indeed, if $\mathfrak{p}$ is of residue class degree $f$ then $q = p^f$ and 
\[ 1 - \alpha_\mathfrak{p}q^{-s} = \underset{j=1}{\overset{f}{\prod}}(1-\zeta_f^j\beta_\mathfrak{p} p^{-s}) \; ,\] where $\zeta_f$ is a primitive $f$-th roof of unity and $\alpha_{\mathfrak{p}} = \beta_\mathfrak{p}^f$. 

Then, the crucial point is that applying the formula to $L(s)$ gives sums of the form 
\[ \underset{m \geqslant 1}{\sum}\underset{j=1}{\overset{f}{\sum}}
(\zeta_f^j\beta_\mathfrak{p})^m \frac{F(m\log p)}{p^m} \log p \; \cdot \]
Using the vanishing of sums of powers of roots of unity, all the terms for $f \nmid m$ disappear, leaving  
\[ \underset{j=1}{\overset{f}{\sum}}\underset{m \geqslant 1}{\sum}
(\zeta_f^j\beta_\mathfrak{p})^{fm} \frac{F(fm\log p)}{p^{fm}} \log p
 =  \underset{m \geqslant 1}{\sum}
\alpha_\mathfrak{p}^{m} \frac{F(m\log q)}{q^m} \log q \; ,\]
and the rest easily follows.
\end{proof}

From this proposition we get the following theorem. Define, for any compact Mestre test function $F(x)$ and real $\lambda >0$,
\[B(K,F,\lambda) = \frac{1}{\delta_K ^2} \exp \left(  M_{\lambda,F} - \frac{2}{n}\underset{\mathfrak{p},m \geqslant 1}{\sum} \left\lfloor 2\sqrt{q^m} \right\rfloor\frac{ | F_\lambda(m\log q) |}{q^{m}}\log q \right) \; ,\] 
where $n$ is the degree of the field $K$ and $\delta_K$ is the root discriminant of $K$. Since $M_{\lambda,F}$ is independent of the field $K$, and by the properties of the root discriminant, this bound is well behaved in towers of number field. 
\begin{theorem}\label{main result} Let $K$ be a number field of degree $n$. Let $E / K$ be an elliptic curve of conductor~$\mathfrak{f}_{E/K}$. For any $\lambda>0$ and $F(x)$ a compact Mestre test function, the following inequality holds: 
\[ B(K,F,\lambda)
\leqslant \left( N_{K/\mathbb{Q}}(\mathfrak{f}_{E/K}) \right) ^{\frac{1}{n}} \cdot \] 
In particular, if the bound $B(K,F,\lambda)$ is greater than $1$, then there are no elliptic curves with everywhere good reduction over $K$.
\end{theorem}
\begin{proof}
We follow the strategy of Mestre used in the previous sections: take $F(x)$ a compact Mestre test function. By Proposition~\ref{equality in number fields} we have
\[ \lambda\underset{\rho}{\sum}\Phi_F(\rho) + 2\underset{\mathfrak{p},m \geqslant 1}{\sum}
b_m(\mathfrak{p})
\frac{F_\lambda(m\log q)}{q^{m}}\log q
+ n M_{\lambda,F} 
=\log(N_{K/\mathbb{Q}}(\mathfrak{f}_{E/K})d_K^2) \cdot  \]
The terms $\Phi_F(\rho)$ are positive since $F(x)$ is a compact Mestre test function. Dividing by $n$, we obtain
\[  \frac{2}{n}\underset{\mathfrak{p},m \geqslant 1}{\sum}
b_m(\mathfrak{p})
\frac{F_\lambda(m\log q)}{q^{m}}\log q
+ M_{\lambda,F} 
\leqslant \log(N_{K/\mathbb{Q}}(\mathfrak{f}_{E/K})\delta_K^2) \cdot  \] The result follows by using the Weil bound to minorate the terms $b_m(\mathfrak{p})$.
\end{proof}

The converse of Theorem~\ref{main result} may not be true: for a given $K$ if
\[  B(K,F,\lambda) \leqslant 1 \] for every compact Mestre test function $F(x)$ and real parameter $\lambda >0$, we cannot deduce the existence of an elliptic curve with everywhere good reduction on $K$. Indeed the existence of an abelian variety of any dimension with everywhere good reduction implies the above bound.

\subsection{Behavior of $B(K,F,\lambda)$ in towers of number field}

To obtain a set $\mathcal{F}_{n-EGR}$ of number fields over which no elliptic curves with everywhere good reduction (and analytic $L$-function) can exist, Theorem~\ref{main result} allows us to simply compute the bound $B(K,F,\lambda)$ for any number fields $K$ and a well chosen set of parameters $F(x)$ and $\lambda$, depending on $K$. If this constant is greater than $1$, we record $K$ in $\mathcal{F}_{n-EGR}$.

This construction of $\mathcal{F}_{n-EGR}$ can be improved in the following way. We observe that if a fields $K$ is in $\mathcal{F}_{n-EGR}$, then its subfields are also in $\mathcal{F}_{n-EGR}$, so we can thake them into account in our computations.  Nevertheless, one could expect, when adding a number field in $\mathcal{F}_{n-EGR}$, that its subfields were already present in the list. Indeed it is unclear that after computing $B(K,F,\lambda)$ for a number field $K$, for some good choice of $\lambda$ and $F(x)$, that there can be a non-trivial extension $L/K$ with $B(L,F,\lambda) > B(K,F,\lambda)$. We prove below that this is possible for fields of class number greater than $1$.


We now examine more closely the behavior of the bound $B(K,F,\lambda)$ under a base extension to~$L/K$.
Fix a prime $\mathfrak{p}$ of $K$ and integer $m \geqslant 1$. We suppose that a Mestre test function $F(x)$ and a parameter $\lambda > 0$ are fixed. When considering a finite extension $L/K$, the term  
\[ \left\lfloor 2\sqrt{(\mathcal{N}\mathfrak{p})^m} \right\rfloor 
\frac{ | F_\lambda(m\log (\mathcal{N}\mathfrak{p})) |}{(\mathcal{N}\mathfrak{p})^{m}}
\log (\mathcal{N}\mathfrak{p})  \] 
is replaced by 
\[ \frac{1}{[L:K]} \underset{\mathfrak{q}\mid \mathfrak{p}}{\sum}
\lfloor 2\sqrt{ (\mathcal{N}\mathfrak{q})^m} \rfloor
\frac{ | F_\lambda(m\log (\mathcal{N}\mathfrak{q})) |}{(\mathcal{N}\mathfrak{q})^{m}}
\log (\mathcal{N}\mathfrak{q}) \]
with $\mathcal{N}\mathfrak{q}$ denoting $N_{L/K}\mathfrak{q}$. First suppose that $\mathfrak{p}$ splits completely in $L$. Then the two previous expressions are equal, since ${N}_{L/K}\mathfrak{q} = {N}_{K/\mathbb{Q}}\mathfrak{p}=q$ for all $\mathfrak{q} \mid \mathfrak{p}$.

Conversely, if $\mathfrak{p}$ is inert in $L$, then the last term is actually 
\[ \left\lfloor 2\sqrt{ (\mathcal{N}\mathfrak{p}) ^{m [L:K]}} \right\rfloor
\frac{ | F_\lambda(m [L:K] \log (\mathcal{N}\mathfrak{p})) |}{(\mathcal{N}\mathfrak{p})^{m [L:K]}} \log (\mathcal{N}\mathfrak{p}) \; , \] 
since ${N}_{L/K}\mathfrak{q} = q^{[L:K]}$. In the case where $|F(x)|$ is a decreasing function, this term is lower than the corresponding term for $K$. This assumption is satisfied when $F(x)$ is the Odlyzko test function.
In general, we have the following:
\begin{lemma}\label{lem for Hilbert extension}Let $K$ be a number field, $L/K$ a finite extension of $K$, $\mathfrak{p}$ a prime ideal of $K$. Let $F(x)$ be a decreasing compact test function and $\lambda>0$. We have:
\[ \left\lfloor 2\sqrt{ (\mathcal{N}\mathfrak{p}) ^m} \right\rfloor
\frac{ | F_\lambda(m\log (\mathcal{N}\mathfrak{p})) |}{(\mathcal{N}\mathfrak{p})^{m}}
\log q  \geqslant 
\frac{1}{[L:K]} \underset{\mathfrak{q}\mid \mathfrak{p}}{\sum}
\left\lfloor 2\sqrt{ \mathcal{N}\mathfrak{q}^m} \right\rfloor
\frac{ | F_\lambda(m\log (\mathcal{N}\mathfrak{q})) |}{(\mathcal{N}\mathfrak{q})^{m}}
\log (\mathcal{N}\mathfrak{q}) \cdot\] 

\end{lemma}

\begin{proof}
Each term of the right is of the form 
\[ f_\mathfrak{q} \left\lfloor 2\sqrt{ (\mathcal{N}\mathfrak{q})^m} \right\rfloor
\frac{ | F_\lambda(m\log \mathcal{N}\mathfrak{q}) |}{(\mathcal{N}\mathfrak{q})^{m}}
\log (\mathcal{N}\mathfrak{p}) \; , \] 
where $f_\mathfrak{q}$ is the inertial degree of $\mathfrak{q}$. 
Since $\mathcal{N}\mathfrak{p} \leqslant \mathcal{N}\mathfrak{q}$, we get 
\[\left\lfloor 2\sqrt{ (\mathcal{N}\mathfrak{q})^m} \right\rfloor
\frac{ | F_\lambda(m\log \mathcal{N}\mathfrak{q}) |}{(\mathcal{N}\mathfrak{q})^{m}}  \leqslant 
\left\lfloor 2\sqrt{ (\mathcal{N}\mathfrak{p})} \right\rfloor
\frac{ | F_\lambda(m\log \mathcal{N}\mathfrak{p}) |}{(\mathcal{N}\mathfrak{p})} \; \cdot\]
The sum of the $f_\mathfrak{q}$ is at most $[L:K]$, hence the result.
\end{proof}
 
Finally we note that the term $\frac{1}{\delta_K ^2}$ is replaced by $\frac{1}{\delta_L ^2}$ when considering the formula over $K$ relative to the formula over an extension $L/K$. These extreme cases should convince the reader that looking for an extension $L/K$ providing a better lower bound on conductors of elliptic curves means looking for an extension where the small primes of $K$ have a high inertia degree, and with a root discriminant close that of $K$. Hence the following theorem:

\begin{theorem}
Let $F$ be a Mestre test function which is decreasing for $x>0$. Fix $\lambda >0$, $K$ a number field and $L$ the Hilbert class field of $K$. Then
\[ B(K,F,\lambda) \leqslant B(L,F,\lambda) \cdot\]
\end{theorem}

\begin{proof}
This simply comes from the fact that $\delta_K=\delta_L$, combined with Lemma \ref{lem for Hilbert extension}.
\end{proof}

\begin{example} Choose $F(x) = (1-|x|)\cos(\pi x) + \frac{1}{\pi}\sin( \pi |x|)$ and let $K = \mathbb{Q}(\sqrt{-23})$. The bound found using this method is $0.45$, with $\lambda =1.84$. This choice of $\lambda$ was obtained testing all $\lambda$ from $1$ to $3$, with an increment of $1/100	$. Taking $L$ the Hilbert class field of $K$ and the same $\lambda$, the lower bound for $L$ is $0.68$, but for $\lambda = 3.58$ we get a bound of $2.01$, proving there can be no elliptic curve with everywhere good reduction over $L$ and therefore over $K$. This reconfirms a result of Cremona \cite[Theorem 5.1]{Cremonagoodreduction}.
\end{example}

\begin{example} Taking the appropriate values for $\lambda$, we obtain that for $f \in \{3, 4, 5, 7, 8, 9, 12, 15 \}$ there is no elliptic curve with everywhere good reduction over $\mathbb{Q}(\zeta_{f})$, computationally confirming the result proved by Schoof \cite[Theorem~1.1]{Schoof2003AbelianVO}. 
\end{example}

If one wants to give statistical results on fields, a good idea is to limit the number of fields we examine by establishing a bound on their root discriminant. We clearly have 
\[ \frac{1}{\delta_K ^2} \exp \left(  M_{\lambda,F} - \frac{2}{n}\underset{\mathfrak{p},m \geqslant 1}{\sum} \left\lfloor 2\sqrt{q^m} \right\rfloor\frac{ | F_\lambda(m\log q) |}{q^{m}}\log q \right)
\leqslant
 \frac{1}{\delta_K ^2} \exp ( M_{\lambda,F} )
 \; ,\] and we note that for all $\lambda >0$, the following inequality holds independently on $F$:

\[M_{\lambda,F} \leqslant 2 \left( \log(2\pi)+ \int_0^{\infty} \frac{e^{-x}}{1-e ^{-x}}-\frac{e^{-x}}{x} \, dx \right)  \; ,\] since $F(x)\leqslant 1$. The above integral equals the Euler's constant $\gamma$, from which we obtain the following bound. 
\begin{proposition} Let $K$ be a number field. If there exists$\lambda >0$ and a compact Mestre test function $F$ such that \[ B(K,F,\lambda) \geqslant 1 \;, \] then 
\[ \delta_K \leqslant 2 \pi e^\gamma \;. \]

\end{proposition}

\subsection{Abelian varieties over number fields}

Following on Mestre's ideas \cite[p~21]{Mestre}, one can use the formulas to obtain results on the reduction of abelian varieties of any dimension. For a number field $K$ of degree $n$ and discriminant $d_K$, and an abelian variety $A/K$ of dimension $g$ over $K$, we attach an $L$-function to the representation on the Tate module $T_\ell(A)$. The Euler factor at $\mathfrak{p}$ is determined by the characteristic conjugacy class of Frobenius at $\mathfrak{p}$, giving the Euler product
\[ L_A(s) = \underset{\mathfrak{p}}{\prod} L_\mathfrak{p}(s) \; ,\] where $ L_\mathfrak{p}(s)$ is a polynomial in $p^{-s}$ of degree $2g$, whose reciprocal roots are of absolute value $\sqrt{ \mathcal N{\mathfrak{p}}}$, in the case of good reduction at $\mathfrak{p}$. In case of mixed bad reduction then the absolute value of the reciprocal roots can also be  $1$ (corresponding to a toric part in the reduction), or even $0$ (that is, the $L$-factor attached to the unipotent part of the reduction is the constant factor $1$). We note $\mathfrak{f}_{A/K}$ the conductor of $A$. The completed $\Lambda$-function is 
\[ \Lambda_A(s) = c^{s/2}(\Gamma_\mathbb{C}(s))^{ng}L(s) \; ,\]
where $ c = N_{K/\mathbb{Q}}(\mathfrak{f}_{A/K})d_K^{2g}$ (see \cite[\S~4]{Serre_facteurs_locaux}). We can now state the following conjecture.
\begin{conjecture} The function $L(s) = L_A\left(s+\frac{1}{2}\right)$ is an analytic $L$-function.
\end{conjecture}

When $\mathfrak{p}$ is of good reduction, the reciprocal roots of the $L_\mathfrak{p}(s)$ can be grouped into $g$ pairs of complex conjugates (see \cite[Proposition 2.1.2]{SerreCurveFiniteFields}). Since they are of absolute value $\sqrt{ \mathcal N{\mathfrak{p}}}$, this leads to a generalization of the inequality $b_m(p) \geqslant - \left\lfloor 2\sqrt{p^m}  \right\rfloor $ of Theorem \ref{Mestre theorem}. We obtain 
\[ b_m(\mathfrak{p}) 
\geqslant 
-  \left\lfloor 2g\sqrt{\mathcal{N}\mathfrak{p}^m}  \right\rfloor  \] 
where $b_m(\mathfrak{p})$ is the term appearing in Theorem~\ref{main result}. By an argument of Serre \cite[Theorem 2.1]{SerreCurveFiniteFields} we have the stronger bound 

\[ b_m(\mathfrak{p})
\geqslant 
- g \left\lfloor 2\sqrt{\mathcal{N}\mathfrak{p}^m}  \right\rfloor .\]
The same inequality holds when $\mathfrak{p}$ is of bad reduction, for the same reason it did in the case of elliptic curves. This, with the form of the functional equation associated to $L_A(s)$, leads to following theorem. The proof is the same as in Theorem \ref{Mestre inequality abelian variety over Q}.

\begin{theorem}\label{Principal th of the article} Let $K$ be a number fields of discriminant $d_K$ and degree $n$. Let $F(x)$ be a even compact Mestre test function and $\lambda >0$. Then 
\[ B(K,F,\lambda)^g
\leqslant \left( N_{K/\mathbb{Q}}(\mathfrak{f}_{A/K}) \right) ^{\frac{1}{n}} \cdot \]
\end{theorem}

Using this theorem in the special case of $K = \mathbb{Q}$, Mestre proved that the conductor of an abelian variety over $\mathbb{Q}$ of dimension $g$, and with analytic $L$-function, is at least $10^g$. We obtain an analogous result in the case of abelian varieties over number fields, hence obtaining a criterion not only for the non-existence of elliptic curves with everywhere good reduction, but for abelian varieties with everywhere good reduction. 

The following table gives the number of fields with bounded root discriminant and computed Mestre bound greater than $1$, sorted by degree.  As we mentioned before, if there is no abelian variety with everywhere good reduction over a number field, then the same is true for any of its subfield. Counting those subfields yields the following table, where the rows corresponding to the degrees $2$ up to $6$ show improved results, displayed in the third column. 

The second column gives the root discriminant bound up to which the computations were completed. When the maximal discriminant considered is smaller than $2 \pi e^\gamma \approx 11.18 $, there might be an exceptional field with higher root discriminant for which the Mestre bound is greater than $1$, but the behavior of the bound makes it highly improbable.  

Finally, between parenthesis we give the number of fields for which our results stay true without supposing GRH. As explained in remark~\ref{without GRH}, we changed the test function we used to get the bound in each computation, obtaining each time a slightly worse bound, hence a lower number of fields with the desired property than when we suppose GRH.

\begin{table}[H]
\begin{tabular}{|c|c|c|c|c|c|c|} \hline

Degree & Root Discriminant bound &  Fields   & Subfields\\ \hline \hline

        2        &  	11.135			&		10   (10)  	& 18 	(18)	\\ \hline

3                &   11.184			&	    6    (6) 	& 11	    	(10)	 \\ \hline

 4               &  11.190       	&		35   (33)	& 48	    (39)  	\\ \hline

  5              &   11.109       	&		40 	 (39) 	& 42	    (39)  	\\ \hline

   6             &   6.995			&		122  (97)   & 126	(97)		\\ \hline

    7            &  8.980  			&		126  (89)  	& 123	(89)		\\ \hline

  8              &   9.202   		&		140  (111)	& 140	(111)		\\ \hline

   9             &   9.953			&		98   (71)	& 98		(71)	\\ \hline

    10           &  8.613  			&		317  (89)   & 317	(89)		\\ \hline

    11           &   11.169     		&		10   (0) 	& 10		(0)	\\ \hline

      12         &   11.186     		&		32  	 (0)		& 32		(0)			\\ \hline

        13       &           		&		0 	(0)		&	0	(0)			\\ \hline

          14     &           		&		0 	(0)		&	0	(0)		\\ \hline

            15   &          			&		0 	(0)		&	0	(0)		\\ \hline

      16         &          			&		0    (0)	   	&	0 	(0)		\\ \hline

\end{tabular}
    \caption{Number of fields with computed bound greater than 1}
    \label{fields without egr}
\end{table}

The detailed analysis of the relation between the bound obtained for a field $K$ and the arithmetic of this field is still being investigated. The full list of fields  with computed bound greater than 1, as well as their LMFDB labels and associated computed bounds, can be found in \cite{DOI_conductor_abelian_varieties}.

\begin{remark}
As mentionned above, it is possible that over a given field $K$, there are no elliptic curves with everywhere good reduction but there exists an abelian variety with everywhere good reduction. Dembelé and Kumar give such examples in~\cite{Dembele_examples_good_reduction}, using the existence of Hilbert modular forms of parallel weight $2$ and level $1$. They only give examples over fields with trivial narrow class number. The same result should hold over $\mathbb{Q }(\sqrt{10})$, which has narrow class number $2$. Indeed, there is no elliptic curve with trivial conductor over $\mathbb{Q }(\sqrt{10})$, but there exists a Hilbert modular form with trivial conductor over $\mathbb{Q }(\sqrt{10})$, and we expect an associated abelian surface with everywhere good reduction to exist. 
\end{remark}

\begin{remark}
Producing conductor bounds in the case of prescribed bad reduction can be treated just as in the previous section in the case of abelian varieties number fields. There are two reasons we decided not to treat this question: the first is the possibility of an elliptic curve having no bad reduction at all on some number field, and we chose to focus on in the paper. The second reason is about the choice of the primes of bad reduction: how to do it for a large set of number fields in an interesting way? For instance: is considering that ``all primes above 2 are supposed to be of bad reduction" is interesting? is it realistic for elliptic curves over number fields of high degree? Or should one try to include randomness in the selection of primes of bad reduction? All these questions are left to future investigation.
\end{remark}

\appendix 

\section{On the choice of good test functions}

In this section, we seek to optimize the choice of a Mestre test function. More explicitly, given a number field $K$ and a parameter $\lambda > 0$, we would like to produce a test function $F(x)$ for which the bound $ B(K,F,\lambda)$ is maximal among all the choices of test functions. We will answer a slightly different problem: firstly, we restrict our attention to test functions in the set 
\[\mathcal{V}= \{g*g \; : \;  g \in \mathcal{C}^2( 
[-1/2,1/2)] , \, g \, \mathrm{even}\} \; \cdot\]
It is easy to verify that any element $F(x)$ in $\mathcal{V}$ such that $F(0)=1$ is a test function. Among this subset of test functions, the Odlyzko test function is optimal in a sense made precise in Poitou~\cite{Poitou}: roughly speaking it is the function with the smallest second derivative at $0$.

We introduce a topology on this space, induced by the norm
\[||F|| = ||F||_\infty + ||F''||_\infty \; \]
on test functions $F(x) \in \mathcal{V}$.
Now, since we fixed the parameter $\lambda$, we will use the coefficient $a_p$ appearing in the explicit formula, and instead of working on $B(K,F,\lambda)$ we focus on the expression 
\[C(K,F,\lambda,(a_\mathfrak{p})_\mathfrak{p})= \frac{2}{n}\underset{\mathfrak{p},m \geqslant 1}{\sum}b_m(\mathfrak{p})\frac{F_\lambda(m\log q)}{q^{m}}\log q+  M_{\lambda,F} \; , \]which is a linear function in $F(x)$, in contrast to $B(K,F,\lambda)$ (note that this is essentially the term appearing in Proposition \ref{equality in number fields}). The dependence in $(a_\mathfrak{p})_\mathfrak{p}$ is contained in the term $b_m(\mathfrak{p})$ in the right hand side. 

To obtain a bound on the conductor of all elliptic curves, it is enough to give a bound on all the $C(K,F,\lambda,(a_\mathfrak{p})_\mathfrak{p})$ for all choices of $(a_\mathfrak{p})_\mathfrak{p}$ with $\mathcal{N}(\mathfrak{p}) \leqslant \exp(\lambda)$, which is a finite set of choices since $a_\mathfrak{p}$ must satisfy the Weil bound $|a_\mathfrak{p}| \leqslant 2 \sqrt{\mathcal{N}_K \mathfrak{p}}$.

We prove now that among all test functions, the ones having a polynomial expression are good enough to maximize $C(K,F,\lambda,(a_\mathfrak{p})_\mathfrak{p})$. We need two lemmas: 

\begin{lemma} Given a field $K$, a finite sequence of coefficients $(a_\mathfrak{p})_\mathfrak{p}$ and a parameter $\lambda >0$, the function $C(K,F,\lambda,(a_\mathfrak{p})_\mathfrak{p})$ is continuous on $\mathcal{V}$ equipped with the above topology.

\begin{proof}

Let $F \in \mathcal{V}$. We have to prove that $C(K,F,\lambda,(a_\mathfrak{p})_\mathfrak{p})$ is bounded if $||F||$ is bounded.
As the parameter $\lambda$ is fixed, the sum 
\[ {\sum}b_m(\mathfrak{p})\frac{F_\lambda(m\log q)}{q^{m}}\log q\] has only finitely many terms and is obviously bounded by $c_1||F||_\infty$ with $c_1$ independent of $F$ (but depending on the field $K$ and the parameter $\lambda$) . We then only need to examine 
\[ M_\lambda(F) = 2 \left( F(0)\log(2\pi)+ \int_0^{\infty} \frac{F_\lambda(x)e^{-x}}{1-e ^{-x}}-F(0)\frac{e^{-x}}{x} \, dx \right) \cdot \]
Without loss of generality we suppose here that $\lambda = 1$. We separate the integral into the integrals from $0$ to $1$ and from $1$ to $\infty$. The second integral
\[ \int_1^{\infty} \frac{F(x)e^{-x}}{1-e ^{-x}}-F(0)\frac{e^{-x}}{x} \, dx\] 
has an absolute value bounded by 
\[ \int_1^{\infty}
 \frac{||F||_\infty e^{-x}}{1-e ^{-x}}
 +||F||_\infty\frac{e^{-x}}{x} \, dx 
:= c_2||F||_\infty \] where $c_2 =  \int_1^{\infty}
 \frac{e^{-x}}{1-e ^{-x}}
 +\frac{e^{-x}}{x} \, dx 
$. For the first integral, we note that
\[ \frac{F(x)e^{-x}}{1-e ^{-x}} - \frac{F(x)}{x} = -F(x) \frac{e^x-x-1}{x(e^x-1)}
 = F(x)(-\frac{1}{2}+O(x))\]
and
\[F(0)\frac{e^{-x}}{x} - F(0)\frac{1}{x} = F(0)(-1+O(x)) \,, \] hence 
\[ 
\int_0^1 \frac{F(x)e^{-x}}{1-e ^{-x}}-F(0)\frac{e^{-x}}{x} \, dx
=
\int_0^1 \frac{F(x)-F(0)}{x} \, dx + \int_0^1 F(x)(-\frac{1}{2}+O(x))+(F(0)(-1+O(x)) \, dx 
\] whose absolute value is smaller than 
\[ \left| \int_0^1 \frac{F(x)-F(0)}{x} \, dx \right| + c_3||F||_\infty \, \cdot\] for some constant $c_3$. 
Since $F(x)$ is even then $F'(0) = 0$, therefore, by repeatedly applying the mean value theorem,
\[ 
\frac{F(x)-F(0)}{x} \leqslant ||F''||_\infty  \; ,
\] 
so that \[ |C(K,F,\lambda,(a_\mathfrak{p})_\mathfrak{p})| \leqslant (c_1+c_2+c_3+2\log(2\pi))||F||_\infty + ||F''||_\infty \, ,  \]

which allows us to establish the continuity of $C(K,F,\lambda,(a_\mathfrak{p})_\mathfrak{p})$ on $\mathcal{V}$.
\end{proof}
\end{lemma}

\begin{lemma} The subset of polynomial functions in $|x|$ forms a dense subspace of $\mathcal{V}$.
\begin{proof} Let $\varepsilon >0$ and let $F(x) \in \mathcal{V}$, say $F(x) = g*g(x)$. Using Weierstrass theorem, we obtain a polynomial function $q(x)$ such that $||q - g'||_{\infty,(-\frac{1}{2},\frac{1}{2})} \leqslant \varepsilon$. One may easily suppose that $q(x)$ is an odd function, as $g'(x)$ is so. By integrating, one gets an even polynomial function $p(x)$ such that $p'(x) = q(x)$ and $||p - q||_{\infty,(-\frac{1}{2},\frac{1}{2})} \leqslant \varepsilon$. 

Now set $\tilde{F}(x) = p*p(x)$. One easily gets 
\[ \mid (F-\tilde{F})(x) \mid  =  \mid (p-q)*(p+q)(x) \mid \leqslant 3 ||q||_\infty \varepsilon \] and similarly, as $F''(x) = (g'*g')(x)$, we have
\[ \mid (F-\tilde{F})''(x) \mid  =  \mid (p'-q')*(p'+q')(x) \mid \leqslant 3 ||q'||_\infty \varepsilon \] if $\varepsilon$ is small enough . This proves the claim, as $p*p(x)$ has a polynomial expression for $x \geqslant 1$.

\end{proof}
\end{lemma}

From these two lemmas, maximizing the operator $C(K,F,\lambda,(a_\mathfrak{p})_\mathfrak{p})$ on $\mathcal{V}$ can be done by using polynomials, since this operator is continuous and polynomials are dense in $\mathcal{V}$. From the above proof it is clear that we can work with functions in $\mathcal{V}$ of the form $p*p$ with $p$ being an even polynomial. As $C(K,F,\lambda,(a_\mathfrak{p})_\mathfrak{p})$ is linear in $F$, thus it is quadratic in $p$.

We illustrate the consequences with an example:  

\begin{example} We look for a function of the form $g * g$ where $g$ is polynomial of the form $g(x) = a_0 + a_2 x^2$ on $[-\frac{1}{2},+\frac{1}{2}]$ that maximizes the conductor lower bound for abelian varieties of dimension $2$ and rank $1$, with multiplicative reduction of dimension $2$ at $p=2$  and additive reduction of dimension $1$ at $p=3$. For the primes of good reduction, we will always suppose that the local factor $L_p(s)$ has the form 
\[ L_p(s) = (1-a_p p^{-s}+p^{1-2s})^{-2} \] 
with $a_p = \lfloor 2 \sqrt{p} \rfloor$. This hypothesis on the form of Euler factors at good primes is heuristic: the explicit formulas suggest that abelian varieties with this kind of Euler factors have the smallest conductor. One may also apply the technique below to all possible forms of the Euler factors at good primes, as there are only a finite number of them for small primes.

A preliminary test using the Odlyzko test function leads to a choice of $\lambda = 2.08$. For this parameter, we then compute the following bounds:
\[ b_{00} = \log(C(K,p_0*p_0,\lambda,(a_p)_{p})) \:,\]
\[ b_{02} = \log(C(K,p_0*p_2,\lambda,(a_{p})_p)) \;,\]
\[ b_{22} = \log(C(K,p_2*p_2,\lambda,(a_p)_p))  \;,\]
where $p_i$ denotes the polynomial function $x^i$ compactly supported on $[-\frac{1}{2},+\frac{1}{2}]$. Then
\[ \log(C(K,g*g,\lambda,(a_p)_{p})) =  b_{00}.a_0^2 +2b_{02}.a_0a_2+b_{22}.a_2^2 \] and we look for the coefficients $a_0$ and $a_2$ that maximize this expression under the condition $g*g(0) =1$. Using the library Gloptypoly of Octave \cite{Octave}, we find 
\[a_0=  1.3780  \quad \mathrm{and} \quad a_2 = -5.656 \; \cdot\]

The following graph shows the optimized function in black and the Odlyzko test function in blue. 

\begin{figure}[!h]
 	\includegraphics[scale=0.4]{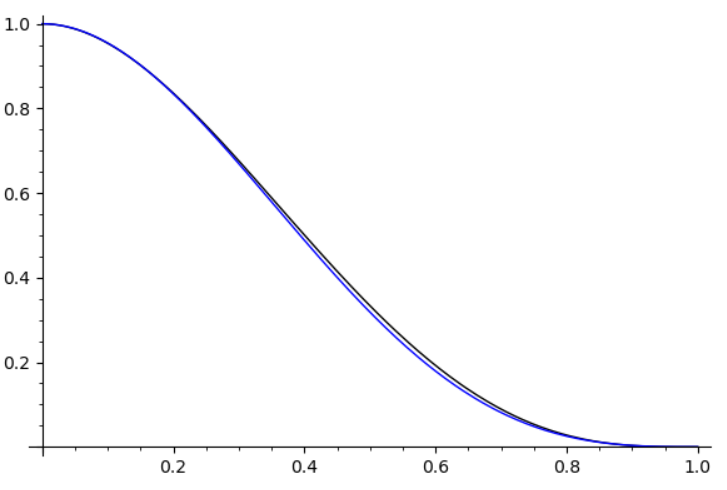} 
\end{figure}
\end{example}

\begin{example} We apply the same method to abelian varieties of dimension $3$, rank $4$, with bad reductions of types $(0,3)$ at $p=2$, type $(2,1)$ at $p=3$ and type $(2,0)$ at $p=5$. We fix $\lambda =2.15$ after some preliminary tests. We obtain $g(x) = a_0+a_2 x^2$ with
\[a_0=1.3846 \quad  \mathrm{and} \quad a_2= -5.791 \]
for which the conductor bound is  $2221.7$. 

The same method applied to the search of a polynomial $g(x) = a_0+a_2 x^2+a_4 x^4$ yields 
\[a_0 = 1.3905, \quad a_2 =  -6.0168 \quad  \mathrm{and}\quad a_4 =  0.99995\] for which the bound is $2225.6$. This can be further improved by increasing the degree of~$g$. 

\end{example}

Remarkably, the Odlyzko test function still yields a better result: under the same conditions the bound for this function is $2232.4$. For this reason, all the computations in the previous sections were done using Odlyzko's function.

\begin{acknowledgments} I would like to thank my advisors Samuele Anni and David Kohel for all the precious remarks they provided, and for the kindness they kept on showing through our discussions.
I heartfully thank Jean-François Mestre, who answered the questions I sent him about his 1986 article which forms the basis of this paper.
I finally give my thanks to Arthur Marmin, without whom numerical optimization of the tests would not have been possible.
\end{acknowledgments}

\bigskip



\bibliographystyle{abbrv}

\end{document}